\documentclass[11pt]{amsart}
\usepackage{amssymb,amsmath,latexsym, enumerate, amscd}
\usepackage[all]{xy}

\theoremstyle{plain}
\newtheorem{theorem}{Theorem}[section]
\newtheorem{corollary}[theorem]{Corollary}
\newtheorem{lemma}[theorem]{Lemma}
\newtheorem{proposition}[theorem]{Proposition}

\theoremstyle{definition}
\newtheorem{definition}[theorem]{Definition}

\newtheorem{remark}[theorem]{Remark}
\newtheorem{Definitions and Notation}[theorem]{Definitions and
Notation}

\numberwithin{equation}{section}

\newcommand{\locoho}[3]{\ensuremath{\operatorname{H}_{#1}^{#2}\left(#3\right)}}

\newcommand{\lra}{\longrightarrow}
\newcommand{\sse}{\subseteq}

\newcommand{\soc}{\operatorname{Soc}}

\newcommand{\seq}[2]{\ensuremath{{#1}_{1}}, \ensuremath{{#1}_{2}}, \ldots, \ensuremath{{#1}_{#2}}{}}

\newcommand{\directlimit}[2]{\ensuremath{\varinjlim_{#1} #2}}

\newcommand{\cln}[3]{\ensuremath{\left(#1:_{#2} #3\right)}}

\newcommand{\Hom}{\operatorname{Hom}}
\newcommand{\Ext}{\operatorname{Ext}}
\newcommand{\KH}{\operatorname{H}}

\newcommand{\limclosure}[1]{\ensuremath{{#1}^{\lim}}}
\newcommand{\hgt}{\ensuremath{\operatorname{ht}}}
\newcommand{\tightclosure}[1]{\ensuremath{{#1}^{*}}}

\begin{document}

\title[Gorenstein Rings]{Gorenstein rings and irreducible parameter
ideals}

\author{Thomas Marley}
\address{Department of Mathematics, University of Nebraska, Lincoln,
NE 68588-0130, USA}
\email{tmarley@math.unl.edu}
\urladdr{http://www.math.unl.edu/\textasciitilde tmarley1}

\author{Mark W. Rogers}
\address{Department of Mathematics, Missouri State University,
Springfield, MO 65897, USA}
\email{markrogers@missouristate.edu}

\author{Hideto Sakurai}
\address{Department of Mathematics, School of Science and
Technology, Meiji University, 214-8571, Japan}
\email{ee78052@math.meiji.ac.jp}

\thanks{The second author was supported for eight weeks during the
summer of 2006 through the University of Nebraska-Lincoln's \it
Mentoring through Critical Transition Points \rm grant (DMS-0354281)
from the National Science Foundation.} 
\date{August 22, 2006}

\begin{abstract}  Given a Noetherian local ring $(R,m)$
it is shown that there exists an integer $\ell$ such that $R$ is
Gorenstein if and only if some system of parameters contained in
$m^{\ell}$ generates an irreducible ideal.  We obtain as a corollary
that $R$ is Gorenstein if and only if every power of the maximal
ideal contains an irreducible parameter ideal.
\end{abstract}

\maketitle

\section{Introduction}

It is well-known that a commutative Noetherian local ring $(R, m)$
is Gorenstein if and only if $R$ is Cohen-Macaulay and some ideal
generated by a system of parameters (called a {\it parameter
ideal}\,) is irreducible. Perhaps less widely known is a result of
Northcott and Rees  which states that if every parameter ideal is
irreducible then $R$ is Cohen-Macaulay \cite[Theorem 1]{NR}.
Hence, $R$ is Gorenstein if and only if every parameter ideal is
irreducible.   There are, however, examples of non-Gorenstein
rings possessing irreducible parameter ideals: $(y)R$ is
irreducible in the local ring $R=\mathbb Q[[x,y]]/(x^2,xy)$, for
example.   In a discussion of the Northcott-Rees result between
the second author and William Heinzer, the following question
arose: {\it If $R$ contains a system of parameters \seq{x}{d} such
that for every positive integer $n$, the ideal $\left(x_{1}^{n},
x_{2}^{n}, \ldots, x_{d}^{n}\right)$ is irreducible, is $R$
necessarily Gorenstein?}

A concept related to this question was studied by Hochster: $R$ is
called \emph{approximately Gorenstein} if every power of $m$
contains an irreducible $m$-primary ideal. While approximately
Gorenstein rings must have positive depth, they need not be
Cohen-Macaulay.  In fact, every complete Noetherian domain is
approximately Gorenstein \cite[Theorem~1.6]{Ho}.   However, our
principal result (Theorem~\ref{main}) shows that if a high enough
power of $m$ contains an irreducible \emph{parameter} ideal then the
ring is Cohen-Macaulay (and hence Gorenstein):

\medskip
\noindent {\bf Theorem:} {\it Let $(R,m)$ be a Noetherian local
ring.  Then there exists an integer $\ell$ such that $R$ is
Gorenstein if and only if some parameter ideal contained in
$m^{\ell}$ is irreducible.}
\medskip

As a consequence, a local ring $R$ is Gorenstein if and only if
every power of the maximal ideal contains an irreducible parameter
ideal, answering the question posed above.  We show that the integer
$\ell$ identified in this theorem may be taken to be the least
integer $\delta = \delta(R)$ such that the canonical map
\[
\Ext^d_R(R/m^{\delta},R) \to \directlimit{}{\Ext^d_R(R/m^n,R)} \cong
H^d_m(R)
\]
is surjective after applying the functor $\Hom_R(R/m,-)$, where
$d=\dim R$.

We note that Theorem~\ref{main} is known in some special cases.
Recent work by Goto-Sakurai \cite{GSa1, GSa2}, Liu-Rogers \cite{LR},
and Rogers \cite{R} has shown that some rings having finite local
cohomologies have eventual constant index of reducibility of
parameter ideals, and these results may be used to prove our Main
Theorem under additional hypotheses. (The index of reducibility of
an ideal $I$ is the number of irreducible ideals appearing in any
irredundant expression of $I$ as an intersection of irreducible
ideals; in the case that $R/I$ has finite length, the index of
reducibility of $I$ is the dimension of $\Hom_R(R/m,R/I)$ as an
$R/m$-vector space, where $m$ denotes the maximal ideal of $R$.) To
be precise, suppose $R$ has finite local cohomologies (that is, the
local cohomology modules \locoho{m}{i}{R} have finite length when $i
\neq d$) and assume that one of the following conditions holds:
Either $R$ is quasi-Buchsbaum (that is, $m \locoho{m}{i}{R} = 0$ for
$i \neq d$), or there is some integer $t$ with $0 < t < d$ such that
$\locoho{m}{i}{R} = 0$ for all $i$ with $i \neq 0, t, d$. Then there
is an integer $\ell$ such that the index of reducibility of every
parameter ideal contained in $m^{\ell}$ is equal to $ \sum_{i =
0}^{d} \binom{d}{i} \dim_{R/m}\Hom_R(R/m,\locoho{m}{i}{R}). $ (This
expression for the eventual constant index of reducibility first
appeared in \cite{GSu}.) Thus, if we further assume that every power
of the maximal ideal contains an irreducible parameter ideal, then
this eventual constant index of reducibility must be $1$: i.e.,
$\sum_{i = 0}^{d} \binom{d}{i} \dim_{R/m} \Hom_R(R/m,
\locoho{m}{i}{R}) = 1$. Since $\locoho{m}{d}{R} \neq 0$, we must
have $\locoho{m}{i}{R} = 0$ for $i < d$. Thus, $R$ is Cohen-Macaulay
and hence Gorenstein.

\section{Main Results} \label{S:MainResults}

As general references for terminology and well-known results, we
refer the reader to \cite{Mat} or \cite{BH}. Throughout, $R$ denotes
a Noetherian ring.  In case $R$ is local with maximal ideal $m$ and
$M$ is an $R$-module, the {\it socle} of $M$ is defined to be $(0:_M
m)=\{x\in M\mid mx=0\}$. The socle of $M$ is denoted by $\soc_R M$,
or simply $\soc M$ if there is no confusion about the ring.  It is
clear that $\soc (-)$ is a left exact covariant functor in a natural
way. We often will identify this functor with the functor
$\Hom_R(R/m,-)$.

\begin{definition}
Let $x_{1}$, \ldots, $x_{r}\in R$ and let $M$ be an $R$-module.
Define
\[
\limclosure{\{ x_{1}, \ldots, x_{r} \}}_M := \bigcup_{n \geq 0}
\cln{\left(x_{1}^{n + 1}, \ldots, x_{r}^{n +
1}\right)M}{M}{x_{1}^{n} \cdots x_{r}^{n}}.
\]
If $M = R$ we write $\limclosure{\{ x_{1}, \ldots, x_{r} \}}$ for
$\limclosure{\{ x_{1}, \ldots, x_{r} \}}_R$.
\end{definition}

We make the following remarks concerning this operation:

\begin{remark} \label{rmk1}
Let $I = (x_{1}, \ldots, x_{r})$ be an ideal of $R$ and let $M$ be
an $R$-module.
\begin{enumerate}[(a)]
\item The set $\limclosure{\{ x_{1}, \ldots, x_{r} \}}_M$ is a
submodule of $M$ containing $IM$.

\item If $\{ x_{1}, \ldots, x_{r} \}$ is a regular sequence on $M$
then $\limclosure{\{ x_{1}, \ldots, x_{r} \}}_M = IM$. (cf.
\cite[Theorem 16.2]{Mat}).

\item Consider the direct system $\{ M/(x_{1}^{n}, \ldots,
x_{r}^{n})M \}_{n \geq 1}$ given by the maps
\[
M/(x_{1}^{s}, \ldots, x_{r}^{s})M \xrightarrow{(x_1\cdots x_r)^{t -
s}} M/(x_{1}^{t}, \ldots, x_{r}^{t})M
\]
for  $1\le s \le t$. Then the kernel of the canonical map
\[
\phi_{t} : M/(x_{1}^{t}, \ldots, x_{r}^{t})M \to
\directlimit{}{M/(x_{1}^{n}, \ldots, x_{r}^{n})M} \cong
\locoho{I}{r}{M}
\]
is $\limclosure{\{ x_{1}^{t}, \ldots, x_{r}^{t} \}}_M/(x_{1}^{t},
\ldots, x_{r}^{t})M$. Hence by (b), if $x_{1}, \ldots, x_{r}$ is a
regular sequence on $M$ then $\phi_{t}$ is injective for every
$t$.

\item If $R$ is local ring of prime characteristic and $x_{1}$,
\ldots, $x_{r}$ and $y_{1}$, \ldots, $y_{r}$ are two minimal
generating sets for $I$ then $\limclosure{\{ x_{1}, \ldots, x_{r}
\}} = \limclosure{\{ y_{1}, \ldots, y_{r} \}}$ \cite[Remark
5.6]{Hu}. In this context, $\limclosure{\{ x_{1}, \ldots, x_{r}
\}}$ is called the \emph{limit closure} of $I$ and is denoted by
\limclosure{I}.

\item If $R$ is a local equidimensional ring of prime
characteristic which is the homomorphic image of a Cohen-Macaulay
ring and $x_{1}, \ldots, x_{r}$ are parameters (i.e., $\hgt
(x_{1}, \ldots, x_{r}) = r$), then $\limclosure{(x_{1}, \ldots,
x_{r})} \sse \tightclosure{I}$, where \tightclosure{I} denotes the
tight closure of $I$ \cite[Theorem 2.3(b)]{Hu}.
\end{enumerate}
\end{remark}

The following proposition is known in a more general setting
\cite[Theorem~5.2.3]{St}, but we include a brief proof for the
reader's convenience.

\begin{proposition} \label{regular}
Let $x_{1}$, \ldots, $x_{r}$ be elements in the Jacobson radical of
$R$ and let $M$ be a finitely generated $R$-module.  The following
conditions are equivalent:
\begin{enumerate}[(a)]
\item $\limclosure{\{ x_{1}, \ldots, x_{r} \}}_M = (x_{1}, \ldots,
x_{r})M$.

\item $x_{1}$, \ldots, $x_{r}$ is a regular sequence on $M$.
\end{enumerate}
\end{proposition}

\begin{proof}
By Remark \ref{rmk1}(b), (b) implies (a).  We prove that (a)
implies (b) by induction on $r$. In the case $r=1$, let $x=x_1$.
Suppose $\limclosure{\{ x \}}_M = (x)M$ and $x\alpha = 0$ for some
$\alpha\in M$. We claim that $\alpha \in (x^{k})M$ for all $k \geq
0$. This is clearly true for $k = 0$, so suppose $\alpha =
x^{k}\beta$ for some $k \geq 0$ and $\beta\in M$. Then
$x^{k+1}\beta = 0$, and thus $\beta\in \limclosure{\{
x\}}_M=(x)M$. Hence, $\alpha \in \left(x^{k + 1}\right)M$. As $x$
is in the Jacobson radical and $M$ is finitely generated,
$\cap_{k} (x^{k})M = 0$ by Krull's Intersection Theorem. Hence,
$\alpha = 0$ and $x$ is a non-zero-divisor on $M$.

Suppose now that $r > 1$.  To complete the proof, we will show the
following:
\begin{enumerate}
\item $\limclosure{\{ x_{1}, \ldots, x_{r - 1} \}}_M = (x_{1},
\ldots, x_{r - 1})M$.

\item $x_{r}$ is a non-zero-divisor on $M/(x_{1}, \ldots, x_{r -
1})M$.
\end{enumerate}
Item (1) will allow us to use the inductive hypothesis to conclude
that $x_{1}$, \ldots, $x_{r - 1}$ is a regular sequence on $M$.

To prove (1), let $\alpha\in \limclosure{\{ x_{1}, \ldots, x_{r -
1} \}}_M$. We claim that for all $k \geq 0$, $\alpha\in (x_{1},
\ldots, x_{r - 1})M + (x_{r}^{k})M$. Again by Krull's Intersection
Theorem, this will imply that $\alpha\in (x_{1}, \ldots, x_{r -
1})M$. The case $k = 0$ is clear, so suppose $\alpha = \omega +
x_{r}^{k}\beta$ where $\omega \in (x_{1}, \ldots, x_{r - 1})M$ and
$\beta\in M$. Thus, $x_{r}^{k}\beta \in \limclosure{\{ x_{1},
\ldots, x_{r - 1} \}}_M$. Hence, there exists $t \geq 0$ such that
\[
(x_{1}  \cdots  x_{r - 1})^{t}x_r^k\beta \in \left(x_{1}^{t + 1},
\ldots, x_{r - 1}^{t + 1}\right)M.
\]
Multiplying by $(x_{1} \cdots x_{r - 1})^{k} x_{r}^{t}$, we obtain
\[
(x_{1} \cdots  x_{r})^{t + k}\beta \in \left(x_{1}^{t + k + 1},
\ldots, x_{r - 1}^{t + k + 1}\right)M \sse \left(x_{1}^{t + k +
1}, \ldots, x_{r}^{t + k + 1}\right)M.
\]
Hence, $\beta \in \limclosure{\{ x_{1}, \ldots, x_{r} \}}_M =
(x_{1}, \ldots, x_{r})M$. Thus, $\alpha \in (x_{1}, \ldots, x_{r -
1})M + (x_{r}^{k + 1})M$.

The proof of (2) is similar:  Suppose $x_{r}\alpha\in (x_{1},
\ldots, x_{r - 1})M$ for some $\alpha\in M$. We claim that $\alpha
\in (x_{1}, \ldots, x_{r - 1})M + (x_{r}^{k})M$ for all $k \geq
0$. Suppose $\alpha = \omega + x_{r}^{k}\beta$ where $\omega \in
(x_{1}, \ldots, x_{r - 1})M$ and $\beta\in M$. Then $x_{r}\alpha =
x_r\omega  + x_{r}^{k + 1}\beta$. Hence, $x_{r}^{k + 1}\beta \in
(x_{1}, \ldots, x_{r - 1})M$. Multiplying by $(x_{1} \cdots x_{r -
1})^{k + 1}$, we obtain that
\[
(x_{1} \cdots  x_{r})^{k + 1}\beta \in (x_{1}^{k + 2}, \ldots,
x_{r - 1}^{k + 2})M \sse (x_{1}^{k + 2}, \ldots, x_{r}^{k + 2})M.
\]
Hence, $\beta \in \limclosure{\{ x_{1}, \ldots, x_{r} \}}_M =
(x_{1}, \ldots, x_{r})M$ and $\alpha \in \left(x_{1}, \ldots, x_{r
- 1}\right)M + \left(x_{r}^{k + 1}\right)M$.
\end{proof}

\begin{corollary} \label{phi}
Let $x_{1}$, \ldots, $x_{r}$ elements of the Jacobson radical of $R$
and let $M$ be a finitely generated $R$-module. For each $t \geq 1$
let $\phi_{t}$ be the canonical map $M/(x_{1}^{t}, \ldots,
x_{r}^{t})M \to \locoho{(x_1,\dots,x_r)}{r}{M}$ as in part (c) of
Remark \ref{rmk1}. The following conditions are equivalent:
\begin{enumerate}[(a)]
\item $\{ x_{1}, \ldots, x_{r} \}$ is a regular sequence on $M$

\item $\phi_{t}$ is injective for some $t \geq 1$.

\item $\phi_{t}$ is injective for all $t \geq 1$.
\end{enumerate}
\end{corollary}

%
%

In the sequel we adopt the following notation:  For a sequence of
elements $\mathbf x=x_1,\dots,x_r$ and $t\in \mathbb N$ we let
$\mathbf x^t$ denote the sequence $x_1^t,\dots,x_r^t$. For $x\in R$
we let $K(x)$ denote the Koszul chain complex $0\to
R\xrightarrow{x}R\to 0$, where the first $R$ is in homological
degree $1$. For the sequence $\mathbf x$ the Koszul chain complex
$K(\mathbf x)$ is defined to be the chain complex $K(x_1)\otimes
\cdots \otimes K(x_{r})$.  For $1 \le s \le t$ there exists a chain
maps $\phi^t_s \colon K(\mathbf x^t)\to K(\mathbf x^s)$ given by
$\phi^t_s = \phi^t_s(x_1)\otimes \cdots \otimes \phi^t_s(x_{r})$,
where $\phi^t_s(x)$ is the chain map
\[
\xymatrix@!C{ K(x^t): \ar[d]^{\phi^t_s(x)} & 0 \ar[r] & R
\ar[r]^{x^{t}}
\ar[d]^{x^{t - s}} & R \ar[r] \ar[d]^{=} & 0\\
K(x^s): & 0 \ar[r] & R \ar[r]^{x^{s}} & R \ar[r] & 0
 }
\]
For an $R$-module $M$, the $i$th Koszul cohomology of $M$ with
respect to $\mathbf x$, denoted $\KH^i(\mathbf x;M)$, is the $i$th
cohomology of $\Hom_R(K(\mathbf x),M)$.  The maps $\phi^t_s$ above
induce chain maps
\[
\Hom_R(K(\mathbf x^s),M) \to \Hom_R(K(\mathbf x^t),M)
\]
for all $1 \le s \le t$.  By \cite[Theorem 2.8]{Gr}, we have
$\directlimit{}{\KH^i(\mathbf x^n;M)}\cong \locoho{(\mathbf
x)}{i}{M}$.

We make the following elementary observations concerning direct
limits:

\begin{remark} \label{rmk2}
Let $\{ M_{n}, \lambda^{n}_{p} \}$ be a direct system of $R$-modules
over a directed index set and let $\phi_{t} : M_{t} \to
\directlimit{}{M_{n}}$ be the canonical maps given by the definition
of the direct limit.
\begin{enumerate}[(a)]
\item If $\directlimit{}{M_{n}}$ is finitely generated then
$\phi_{t}$ is surjective for all sufficiently large $t$.

\item If $A$ is a finitely presented $R$-module then
$\Hom_R(A,\directlimit{}{M_n})\cong \directlimit{}{\Hom_R(A,M_n)}$.
\end{enumerate}
\end{remark}

\begin{proof}  Part (a) is an easy consequence of
\cite[Theorem 2.17]{Ro}.  For part (b), see Exercise 26, Chapter III
of \cite{L}.
\end{proof}

\begin{definition} \label{R:ell} Let $(R,m)$ be a local ring,
let $M$ be a finitely generated $R$-module, and let $i\ge 0$. By
applying $\Ext^{i}_{R}(-, M)$ to the system of surjections
\[
\cdots \to R/m^{3} \to R/m^{2} \to R/m
\]
we obtain a direct system whose limit is
$\directlimit{}{\Ext^i_R(R/m^n,M)} \cong \locoho{m}{i}{M}$. By
Remark~\ref{rmk2}(b),
\[
\directlimit{}{\soc \Ext^i_R(R/m^n,M)}\cong\soc
\directlimit{}{\Ext^i_R(R/m^n,M)}\cong \soc \locoho{m}{i}{M}.
\]
Since $\locoho{m}{i}{M}$ is Artinian, $\soc \locoho{m}{i}{M}$ is
finitely generated. Hence, by Remark~\ref{rmk2}(a) there exists a
smallest nonnegative integer $\ell_{i}(M)$ such that the map $\soc
\Ext^i_R(R/m^t,M)\to \soc \locoho{m}{i}{M}$ is surjective for all
$t\ge \ell_i(M)$.
\end{definition}

The following proposition is essentially \cite[Lemma 3.12]{GSa1}. A
complete proof is given in Section~\ref{GSProof}.

\begin{proposition} \cite[Lemma 3.12]{GSa1} \label{GotoSak}
Let $(R,m)$ be a Noetherian local ring and let $M$ be a finitely
generated $R$-module. For $i\ge 0$ and all $m$-primary ideals
$q=(x_1,\dots,x_r)=(\mathbf x)$ contained in $m^{\ell_i(M)}$ the map
\[
\soc \KH^i(\mathbf x;M)\to \soc \locoho{m}{i}{M}
\]
induced by the canonical map $\KH^i(\mathbf x;M)\to
\directlimit{}{\KH^i(\mathbf x^n;M)}$ is surjective.
\end{proposition}

We now proceed with the proof of our main result:

\begin{theorem} \label{main}
Let $(R, m)$ be a Noetherian local ring of dimension $d$ and let
$\ell=\ell_d(R)$. Then $R$ is Gorenstein if and only if some
parameter ideal contained in $m^{\ell}$ is irreducible.
\end{theorem}

\begin{proof}
It suffices to show that if there exists a system of parameters
$\mathbf x=x_1,\dots,x_d$ contained in $m^{\ell}$ which generates an
irreducible ideal then $R$ is Cohen-Macaulay (and hence Gorenstein).
Let $\phi = \phi_{1}$ denote the canonical homomorphism from
$H^d(\mathbf x; R)\cong R/(\mathbf x)$ to
$\directlimit{}{H^d(\mathbf x^n;R)}\cong \locoho{m}{d}{R}$. By
Remark \ref{rmk1}(c) we have an exact sequence
\[
0 \lra \frac{\limclosure{\{\mathbf x \}}}{(\mathbf x)} \lra
\frac{R}{(\mathbf x)} \overset{\phi}{\lra} \locoho{m}{d}{R}.
\]
Applying the socle functor and using Proposition~\ref{GotoSak} we
obtain the exact sequence
\[
0 \lra \soc{\frac{\limclosure{\{\mathbf x \}}}{(\mathbf x)}} \lra
\soc{\frac{R}{(\mathbf x)}} \lra \soc{\locoho{m}{d}{R}} \lra 0.
\]
Since \locoho{m}{d}{R} is a nonzero Artinian module, it has a
nonzero socle. Since $(\mathbf x)$ is irreducible, $R/(\mathbf x)$
has a one-dimensional socle. Hence, $\soc (\limclosure{\{ \mathbf x
\}}/(\mathbf x))=0$, which implies $\limclosure{\{ \mathbf x
\}}=(\mathbf x)$. By Proposition \ref{regular}, we see that $\mathbf
x$ is a regular and hence $R$ is Cohen-Macaulay.
\end{proof}

\begin{corollary} \label{maincor}
Let $(R,m)$ be a Noetherian local ring.  Then $R$ is Gorenstein if
and only if every power of the maximal ideal contains an irreducible
parameter ideal.
\end{corollary}

\begin{proof}
Immediate from Theorem \ref{main}.
\end{proof}

\section{A Proof of Proposition~\ref{GotoSak}} \label{GSProof}

The proof of this proposition as given in \cite{GSa1}, while
illuminating, is quite terse. Since this result is crucial to our
paper, and indeed crucial for all recent research on the index of
reducibility of parameter ideals, we give a more detailed proof in
this section. Throughout, $R$ denotes a Noetherian ring. We begin
with a lemma:
%

\begin{lemma} \label{resolution} Let $\mathbf x = x_{1}, \ldots,
x_{r}$ be a sequence of elements from $R$ and let $I = (\mathbf
x)R$. Then there exist a family of complexes $\{F(t)\}_{t \ge 1}$
and chain maps $\alpha(t) \colon K(\mathbf x^t)\to F(t)$ and
$\beta(t + 1) \colon F(t + 1) \to F(t)$ such that for each $t \ge 1$
\begin{enumerate}
\item $F(t)$ is a free resolution of $R/I^t$ and each $F(t)_i$ is
finitely generated;

\item $F(t)_0=R$;

\item $\alpha(t)_0$ and $\beta(t + 1)_0$ are the identity maps;

\item the diagram
\[
\begin{CD}
K(\mathbf x^{t + 1}) @>\phi(t + 1)>> K(\mathbf x^t) \\
@VV\alpha(t + 1)V @VV\alpha(t)V \\
F(t + 1) @>\beta(t + 1)>> F(t)
\end{CD}
\]
commutes, where $\phi(t + 1)$ is the chain map $\phi^{t + 1}_t$
defined in Section~\ref{S:MainResults}.
\end{enumerate}
\end{lemma}

\begin{proof}  We use induction on $t$.  Choose $F(1)$ to be any
minimal free resolution of $R/I$ and $\alpha(1) \colon K(\mathbf x)
\to F(1)$ any lifting of $\text{id}_{R/I} \colon H_0(\mathbf x) \to
H_0(F(1))$.  Suppose $t \ge 1$ and for all $1 \le k \le t$ there
exists resolutions $F(k)$ and chain maps $\alpha(k)$ and $\beta(k)$
which have the desired properties. We will construct $F(t + 1)$,
$\alpha(t + 1)$ and $\beta(t + 1)$.   First, we simplify notation:
Let $G := F(t)$, $C := K(\mathbf x^{t + 1})$, and $\gamma =
\alpha(t) \phi(t + 1)$.  We need to construct a resolution $F = F(t
+ 1)$ of $R/I^t$ and chain maps $\alpha = \alpha(t + 1) \colon C \to
F$ and $\beta = \beta(t + 1) \colon F \to G$ such that $\gamma =
\alpha \beta$. The proof of this is a variation on the Horseshoe
Lemma \cite[Lemma 6.20]{Ro}.  Let $F_0=R$ and $\beta_0 = \alpha_0 =
\operatorname{id}_R$. Suppose for some $k \ge 0$ there exists a
commutative diagram of the form
\[
\begin{CD}
C_{k + 1} @>\partial_{k + 1}>> C_k @>\partial_k>> C_{k - 1} @>>>
\cdots @>\partial_1>> C_0 @>\partial_0>> R/(\mathbf x^{t + 1}) @>>>
0
\\
@. @VV\alpha_kV @VV\alpha_{k - 1}V @. @V\alpha_0VV @V\alpha_{-1}VV
@. \\
@. F_k @>\partial'_k>> F_{k - 1} @>>> \cdots @>\partial'_1>> F_0
@>\partial'_0>> R/I^{t + 1} @>>> 0 \\
@. @VV\beta_kV @VV\beta_{k - 1}V @. @V\beta_0VV @V\beta_{-1}VV @. \\
G_{k + 1} @>\partial''_{k + 1}>> G_k @>\partial''_k>> G_{k - 1}
@>>> \cdots @>\partial''_1>> G_0 @>\partial''_0>> R/I^t @>>> 0
\end{CD}
\]
where the middle row is exact and  $F_i$ is a finitely generated
free module and $\gamma_i = \alpha_i \beta_i$ for all $i \le k$. (In
the diagram, $\alpha_{-1}$ and $\beta_{-1}$ denote the natural
surjections.) Let $u_1, \ldots, u_s \in F_k$ be generators for $\ker
\partial'_k$ and $w_1, \ldots, w_z \in G_{k + 1}$ be generators for
$\ker \partial''_{k + 1}$.  Let $F_{k + 1}$ be a free $R$-module of
rank $s + z$ and $a_1, \ldots, a_s, b_1, \ldots, b_z$ a basis for
$F_{k + 1}$. Define $\partial'_{k + 1} \colon F_{k + 1} \to F_k$ by
$\partial'_{k + 1}(a_i) = u_i$ for $1 \le i \le s$ and $\partial'_{k
+ 1}(b_i) = 0$ for $1 \le i \le z$. Clearly, $\operatorname{im}
\partial'_{k + 1} = \ker \partial'_k$.  Choose $c_1, \ldots, c_s \in
G_{k + 1}$ such that $\partial''_{k + 1}(c_i) = \beta_k(u_i) =
\beta_k(\partial'_{k + 1}(a_i))$ for $1 \le i \le s$. Define
$\beta_{k + 1} \colon F_{k + 1} \to G_{k + 1}$ by $\beta_{k +
1}(a_i) = c_i$ for $1 \le i \le s$ and $\beta_{k + 1}(b_i) = w_i$
for $1 \le i \le z$. Evidently, $\partial''_{k + 1} \beta_{k + 1} =
\beta_k \partial'_{k + 1}$. Now let $e_1, \ldots, e_p$ be a basis
for $C_{k + 1}$.  Choose $f_1, \ldots, f_p \in F_{k + 1}$ such that
$\partial'_{k + 1}(f_i) = \alpha_k \partial_{k + 1}(e_i)$ for $1 \le
i \le p$.  Then for $1 \le i \le p$,
\begin{align}
\partial''_{k + 1} \gamma_{k + 1}(e_i) &= \gamma_{k} \partial_{k +
1}(e_i) \notag\\
&= \beta_k \alpha_k \partial_{k + 1}(e_i) \notag\\
&= \beta_k \partial'_{k + 1}(f_i) \notag\\
&= \partial''_{k + 1} \beta_{k + 1}(f_i). \notag
\end{align}
Hence, $\gamma_{k + 1}(e_i) - \beta_{k + 1}(f_i) \in \ker
\partial''_{k + 1}$ for $1 \le i \le p$.  For $1 \le i \le p$ let
$v_i \in R b_1 + \cdots + R b_z$ be such that $\beta_{k + 1}(v_i)
= \gamma_{k + 1}(e_i) - \beta_{k + 1}(f_i)$.  Finally, define
$\alpha_{k + 1} \colon C_{k + 1} \to F_{k + 1}$ by $\alpha_{k +
1}(e_i) = f_i + v_i$ for $1 \le i \le p$.  It is easily verified
that $\partial'_{k + 1} \alpha_{k + 1} = \alpha_k \partial_{k +
1}$ and $\gamma_{k + 1} = \beta_{k + 1} \alpha_{k + 1}$.
\end{proof}

\smallskip

Let $N$ be an arbitrary $R$-module, let $\mathbf x=x_1,\dots,x_r$,
and let $I = (\mathbf x)$ as in Lemma~\ref{resolution}.  Applying
$\Hom_R(-,N)$ to the commutative diagram in part (4) of this lemma,
we get for all $t\ge 1$ a commutative square of cochain complexes
\[
\begin{CD}
\Hom_R(F(t),N) @>\beta(t+1,N) >> \Hom_R(F(t+1),N) \\
@V\alpha(t,N) VV @V \alpha(t+1,N) VV \\
\Hom_R(K(\mathbf x^t),N) @>\phi(t+1,N) >> \Hom_R(K(\mathbf
x^{t+1}),N).
\end{CD}
\]
Taking $i$th cohomologies, we have for all $t\ge 1$ the commutative
diagram
\begin{equation} \label{D:ExtKosLoc}
\begin{CD}
\Ext^{i}_{R}(R/I^t, N) @>\beta^i_{t + 1}(N)>> \Ext^{i}_{R}(R/I^{t
+
1}, N) @>>> \directlimit{}{\Ext^{i}_{R}(R/I^n, N)}\\
@V\alpha^{i}_{t}(N)VV @V\alpha^{i}_{t + 1}(N)VV
@V\directlimit{}{\alpha^{i}_{n}(N)}VV \\
\KH^{i}(\mathbf x^{t}; N) @>\phi^{i}_{t+1}(N)>> \KH^{i}(\mathbf x^{t
+ 1}; N) @>>> \directlimit{}{\KH^{i}(\mathbf x^{n}; N)}.
\end{CD}
\end{equation}

\begin{lemma} \label{iso}  For all $i\ge 0$ the map $\directlimit{}{\alpha^i_n(N)}$
of Diagram~\eqref{D:ExtKosLoc} is an isomorphism.
\end{lemma}

\begin{proof} For $i = 0$, we have the diagram
\begin{equation*}
\begin{CD}
\Hom_{R}(R/I^t, N) @>\beta^0_{t + 1}(N)>> \Hom_{R}(R/I^{t +
1}, N)\\
@V\alpha^{0}_{t}(N)VV @V\alpha^{0}_{t + 1}(N)VV\\
\KH^{0}(\mathbf x^{t}; N) @>\phi^{0}_{t + 1}(N)>> \KH^{0}(\mathbf
x^{t + 1}; N).
\end{CD}
\end{equation*}
In this case, $\alpha^0_t(N)$ is injective for all $t$, so
\directlimit{}{\alpha^{0}_{n}(N)} is injective, and since
$\{I^n\}$ and $\{(\mathbf x^n)\}$ are cofinal, the map
\directlimit{}{\alpha^{0}_{n}(N)} is surjective.

Suppose $i\ge 1$ and $\directlimit{}{\alpha^j_n}(M)$ is an
isomorphism for all $0\le j\le i-1$ and all $R$-modules $M$. Let $E$
be the injective hull of $N$ and consider the short exact sequence
\[
0\to N\to E\to C\to 0
\]
where $C=E/N$.  Since $F(t)$ and $K(\mathbf x^t)$ are complexes of
free modules, we have for all $t\ge 1$ a commutative diagram
\[
\xymatrix{
0 \ar[r] & \Hom_R(F(t),N) \ar[r] \ar[d]_{\alpha(t,N)} &
\Hom_R(F(t),E) \ar[r] \ar[d]_{\alpha(t,E)} & \Hom_R(F(t),C)
\ar[r] \ar[d]_{\alpha(t,C)} & 0 \\
0 \ar[r] & \Hom_R(K(\mathbf x^t),N) \ar[r] & \Hom_R(K(\mathbf
x^t),E) \ar[r] & \Hom_R(K(\mathbf x^t),C) \ar[r] & 0 }
\]
where the rows are short exact sequences of cochain complexes.  We
also note that $\beta(t+1)$ and $\phi(t+1)$ induce maps  from this
diagram to the same diagram but with $t$ replaced everywhere by
$t+1$.  (This would be represented by a 3-dimensional commutative
diagram.) Applying the long exact sequence on cohomology to this
diagram, we obtain for $t\ge 1$ a commutative diagram with exact
rows
\[
\begin{CD}
\Ext^{i-1}_R(R/I^t,E) @>>> \Ext^{i-1}_R(R/I^t,C) @>>>
\Ext^i_R(R/I^t,N) @>>> 0 \\
@V \alpha^{i-1}_t(E) VV @V \alpha^{i-1}_t(C) VV @V \alpha^{i}_t(N)
VV
@V\alpha^i_t(E) VV \\
\KH^{i-1}(\mathbf x^t;E) @>>> \KH^{i-1}(\mathbf x^t;C) @>>>
\KH^i(\mathbf x^t;N) @>>> \KH^i(\mathbf x^t;E).
\end{CD}
\]
Taking direct limits of these diagrams and using that
$\directlimit{}{\KH^i(\mathbf x^n;E)}=0$ for all $i\ge 1$
(\cite[Proposition 2.6]{Gr}), we obtain a commutative diagram with
exact rows
\[
\begin{CD} \directlimit{}{\Ext^{i-1}_R(R/I^n,E)} @>>>
\directlimit{}{\Ext^{i-1}_R(R/I^n,C)} @>>>
\directlimit{}{\Ext^i_R(R/I^n,N)} @>>> 0 \\
@V \directlimit{}{\alpha^{i - 1}_n(E)} VV @V
\directlimit{}{\alpha^{i-1}_n(C)} VV  @V
\directlimit{}{\alpha^{i}_n(N)} VV @. \\
\directlimit{}{\KH^{i-1}(\mathbf x^n;E)} @>>>
\directlimit{}{\KH^{i-1}(\mathbf x^n;C)} @>>>
\directlimit{}{\KH^i(\mathbf x^n;N)} @>>> 0.
\end{CD}
\]
Since $\directlimit{}{\alpha^{i-1}_n(E)}$ and
$\directlimit{}{\alpha^{i-1}_n(C)}$ are isomorphisms (by the
induction hypothesis), we conclude that
$\directlimit{}{\alpha^{i}_n(N)}$ is an isomorphism.
\end{proof}

We now proceed with:

\medskip

{\it Proof of Proposition~\ref{GotoSak}:}  With the notation as in
the statement of the Proposition, fix $i\ge 0$ and let
$\ell=\ell_i(M)$. Let $q=(\mathbf x)$ be a parameter ideal contained
in $m^{\ell}$. By applying the socle functor to
Diagram~\eqref{D:ExtKosLoc} (with $q$ in place of $I$) we obtain the
commutative diagram
\begin{equation} \label{D:SocExtKosLoc}
\begin{CD}
\soc \Ext^i_R(R/q,M) @>>>
\soc \directlimit{}{\Ext^i_R(R/q^n,M)}\\
@V\soc \alpha^i_1(M)VV @V\soc \directlimit{}{\alpha^i_n(M)}VV\\
\soc \KH^i(\mathbf x;M) @>>> \soc \directlimit{}{\KH^i(\mathbf
x^n;M)}.
\end{CD}
\end{equation}
By Lemma~\ref{iso}, the right arrow is an isomorphism.  Let $J$ be
an injective resolution of $M$. For each $t \ge 1$ the commutative
square of natural surjections
\begin{equation} \label{D:QsAndMs}
\begin{CD}
R/q^{t+1} @>>> R/m^{\ell (t+1)} \\
@VVV @VVV \\
R/q^{t} @>>> R/m^{\ell t}
\end{CD}
\end{equation}
induces the commutative diagram of cochain complexes
\begin{equation} \label{D:QsAndMsWithHom}
\begin{CD}
\Hom_R(R/m^{\ell t}, J) @>\varphi_t >> \Hom_R(R/q^t, J) \\
@VVV @VVV \\
\Hom_R(R/m^{\ell (t+1)}, J) @>\varphi_{t + 1}>> \Hom_R(R/q^{t+1}, J)
\\
@VVV @VVV \\
\directlimit{}{\Hom_R(R/m^{\ell n}, J)}
@>\directlimit{}{\varphi_n}>> \directlimit{}{\Hom_R(R/q^n,J)}.
\end{CD}
\end{equation}
The map \directlimit{}{\varphi_n} is an isomorphism.  Indeed, since
the maps in Diagram~\eqref{D:QsAndMs} are surjections, the maps in
the top square of Diagram~\eqref{D:QsAndMsWithHom} are injections,
and thus \directlimit{}{\varphi_n} is an injection. Also, since $q$
is $m$-primary, $\{q^n\}$ and $\{m^{\ell n}\}$ are cofinal. Hence,
\directlimit{}{\varphi_n} is surjective.

Taking $i$th cohomologies, applying the socle functor, and using
Remark \ref{rmk2}(b), we obtain the commutative diagram
\begin{equation*}
\begin{CD}
\soc \Ext^i_R(R/m^{\ell},M) @>>> \soc \Ext^i_R(R/q,M) \\
@VVV @VVV \\
\directlimit{}{\soc \Ext^i_R(R/m^{\ell n},M)} @>>>
\directlimit{}{\soc \Ext^i_R(R/q^n,M)}.
\end{CD}
\end{equation*}
Now, the left vertical arrow is surjective by the definition of
$\ell$ and \cite[Exercise~2.43]{Ro}.  As the bottom arrow is an
isomorphism, we conclude that the right arrow is also surjective.
We complete the proof of the proposition by noting that since the
top and right maps in Diagram~\eqref{D:SocExtKosLoc} are
surjective, so is the bottom map. \qed

\bigskip
\noindent {\bf Acknowledgment:} The authors would like to thank
Craig Huneke for offering an important insight to the proof of
Proposition \ref{regular}, and also William Heinzer and Jung-Chen
Liu for careful readings of this manuscript.


\begin{thebibliography}{Mat}
\bibitem[BH]{BH}
Bruns,~W. and Herzog,~J., {\em Cohen-Macaulay Rings}, Cambridge
Studies in Advanced Mathematics no. {\bf 39}, Cambridge, Cambridge
University Press, 1993.

\bibitem[GSa1]{GSa1}
Goto,~S. and Sakurai,~H., {\em The equality $I^2 = QI$ in
Buchsbaum rings}, Rend. Sem. Mat. Univ. Padova \textbf{110}
(2003), 25--56.

\bibitem[GSa2]{GSa2}
Goto,~S. and Sakurai,~H., \emph{Index of Reducibility of Parameter
Ideals for Modules Possessing Finite Local Cohomology Modules},
Preprint.

\bibitem[GSu]{GSu}
Goto,~S. and Suzuki,~N., {\em Index of Reducibility of Parameter
Ideals in a Local Ring}, J. Alg. \textbf{87} (1984), 53--88.

\bibitem[Gr]{Gr} Grothendieck,~A., {\em Local Cohomology}, notes
by R. Hartshorne, Lect. Notes Math. no. 41, Springer, Berlin,
1966.

\bibitem[HH]{HH} Hochster,~M. and Huneke,~C., {\em Tight closures of
parameter ideals and splitting in module-finite extensions}, J. Alg.
Geom. {\bf 3} (1994), 599--670.

\bibitem[Ho]{Ho} Hochster,~M., {\em Cyclic purity versus purity in
excellent Noetherian rings}, Trans. Am. Math. Soc. {\bf 231} (1977),
no. 2, 463--488.

\bibitem[Hu]{Hu} Huneke,~C., {\em Tight closure, parameter ideals,
and geometry}, in Six Lectures on Commutative Algebra (J. Elias,
J.M. Giral, R.M. Mir\'o-Roig, and S. Zarzuela, eds), Progress in
Mathematics, vol. 166, Birkh\"auser Verlag, Basel, 1998, 187--239.

\bibitem[L]{L} Lang,~S., {\em Algebra}, 3rd ed. (revised),
Graduate Texts in Mathematics no. 211, Springer, New York, 2002.

\bibitem[LR]{LR} Liu,~J.~C. and Rogers,~M., {\em The index of
reducibility of parameter ideals and mostly zero finite local
cohomologies}, Comm. Alg., to appear.

\bibitem[Mat]{Mat}
Matsumura,~H., {\em Commutative Ring Theory}, Cambridge Studies in
Advanced Mathematics no. {\bf 8}, Cambridge, Cambridge University
Press, 1986.

\bibitem[NR]{NR}
Northcott,~D.~G. and Rees,~D., {\em Principal Systems}, Quart. J.
Math. \textbf{8} (1957), 119-127.

\bibitem[R]{R}
Rogers,~M., {\em The index of reducibility for parameter ideals in
low dimension}, J. Alg., \textbf{278/2} (2004), 571--584.

\bibitem[Rot]{Ro} Rotman,~J., {\em An Introduction to Homological
Algebra}, Academic Press, Orlando, FL, 1979.

\bibitem[St]{St} Strooker,~J.~R., \emph{Homological questions in
local algebra}, London Math. Soc. Lecture Note series \textbf{145},
Cambridge Univ. Press, Cambridge 1990.
\end{thebibliography}
\end{document}